\documentclass[journal]{IEEEtran}

\usepackage{amsmath}
\usepackage{amssymb}
\usepackage{graphicx}
\usepackage{bm}
\usepackage{placeins}
\usepackage{siunitx}

\begin{document}

\title{Advanced Scenario Creation Strategies for Stochastic Economic Dispatch with Renewables}

\author{Ryan~N.~King, Matthew~Reynolds, Devon Sigler, and~Wesley~Jones
\thanks{R. N. King, M. Reynolds, D. Sigler, and W. Jones are with the National Renewable Energy Laboratory (NREL), Golden, CO, 80401 USA, e-mail: ({ryan.king, matthew.reynolds, devon.sigler, wesley.jones}@nrel.gov).}
\thanks{Manuscript received June 6, 2018.}}

\markboth{IEEE Transactions on Power Systems, 2018}
{Shell \MakeLowercase{\textit{et al.}}: Bare Demo of IEEEtran.cls for IEEE Journals}

\maketitle

\begin{abstract}
Real-time dispatch practices for operating the electric grid in an economic and reliable manner are evolving to accommodate higher levels of renewable energy generation.  In particular, stochastic optimization is receiving increased attention as a technique for handling the inherent uncertainty in wind and solar generation.  The typical two-stage stochastic optimization formulation relies on a sample average approximation with scenarios representing errors in forecasting renewable energy ramp events.  Standard Monte Carlo sampling approaches can result in prohibitively high-dimensional systems for optimization, as well as a poor representation of extreme events that challenge grid reliability.  We propose two alternative scenario creation strategies, importance sampling and Bayesian quadrature, that can reduce the estimator's variance.  Their performance is assessed on a week's worth of 5 minute stochastic economic dispatch decisions for realistic wind and electrical system data.  Both strategies yield more economic solutions and improved reliability compared to Monte Carlo sampling, with Bayesian quadrature being less computationally intensive than importance sampling and more economic when considering at least 20 scenarios.
\end{abstract}


\section*{Nomenclature}
\addcontentsline{toc}{section}{Nomenclature}
\begin{IEEEdescription}[\IEEEusemathlabelsep\IEEEsetlabelwidth{$V_1,V_2,V_3$}]
\item[$x_{g}$]Dispatch of thermal generators
\item[$c_g$]Cost of thermal generator $g$
\item[$L\left(\cdot\right)$]Loss function
\item[$\omega_{w}$]Dispatch of wind generator
\item[$c_w$]Cost of wind plant $w$
\item[$\omega_{w}^{fcst}$]Forecasted wind generation
\item[$\xi_{w}$]Error in wind forecast
\item[$x_{g}^{max}$]Max output of generator $g$
\item[$x_{g}^{min}$]Min output of generator $g$
\item[$R_g^{up}$]Max ramp up of generator $g$
\item[$R_g^{down}$]Max ramp down of generator $g$
\item[$G$]Set of thermal generators
\item[$T$]Set of timesteps
\item[$W$]Set of wind plants
\item[$S$]Set of wind forecast error scenarios
\item[$D$]Set of load buses
\item[$d_{q}$]Value of loads
\item[$y^+_{q}$]Loss-of-load
\item[$y^-_{q}$]Excess capacity
\item[$c^+_{q}$]Cost of Loss-of-load
\item[$c^-_{q}$]Cost of excess capacity
\item[$\omega^{spl}_{w}$]Spilled wind
\item[$c^{spl}_{w}$]Cost of spilled wind
\end{IEEEdescription}

\section{Introduction}
\IEEEPARstart{M}{odern} practices for economic dispatch of the electrical grid are rooted in deterministic problem formulations.  The increasing deployment of stochastic renewable energy generators, such as wind and solar, requires dispatch techniques that are consistent with the reliability and economic operating requirements of the electric grid.  Stochastic optimization techniques have received increased attention from the power systems community, and are a promising way to address the randomness in renewable energy generation.  This paper explores two novel approaches to scenario-based stochastic optimization that yield lower costs and reduced reserve requirements, facilitating the grid integration and further deployment of renewable energy generators.

As access to high performance computing resources grows, grid operators are able to consider increasingly sophisticated optimization problems in their dispatch decisions.  These complexities can include considering AC optimal power flow (ACOPF), larger grids with thousands of buses, multiple stages or time periods, and stochastic optimization techniques.  These techniques can improve grid reliability in the face of uncertain generation, and help remove obstacles to further renewable energy deployment.

An emerging approach to dispatching generators in the face of uncertain renewable generation is to formulate and solve a two-stage stochastic programming problem with recourse.  In two-stage recourse problems a first stage decision is made, followed by a recourse stage where adjustments can be made after a random variable is realized.  In the context of economic dispatch for the grid, this involves making initial generator dispatch decisions that minimize the combined costs of the first stage and the expectation of the recourse costs.  The expected recourse costs are calculated with respect to the distribution of the uncertain renewable generation, and are typically estimated with a sample average approximation using Monte Carlo techniques \cite{constantinescu2011computational,soroudi2013decision}.  However, Monte Carlo techniques exhibit slow convergence with respect to the number of samples and draw samples that are not necessarily a good test of grid robustness in the event of large ramps in generation, creating a need for improved sampling techniques.

In this article we report on two techniques for creating scenarios of errors in forecasting wind generation that improve the convergence rates and accuracy when calculating the expected recourse costs, ultimately achieving more economic dispatch solutions.  The created scenarios also contain ramp events that are more challenging for the grid, resulting in dispatch decisions that produce more reliable grid operation.  The cost savings and reliability improvements are demonstrated in a week-long simulation of 5 minute dispatch decisions on realistic wind and grid data.

In Section II we provide a background on stochastic programming techniques for two-stage recourse problems, the standard sample average approximation, and power grid applications in economic dispatch.  In Section III we introduce new scenario creation strategies using importance sampling and Bayesian quadrature.  In Section IV we demonstrate the implementation of these techniques in a stochastic economic dispatch framework and examine their performance.  Finally, we conclude in Section V with a discussion of the results and extensions to more complex multi-stage problems and more complex network topologies.

\section{Two-Stage Stochastic Programming}
\subsection{General Problem Formulation}
The general two-stage stochastic programming problem \cite{shapiro2009lectures,birge2011introduction} is defined as:
\begin{equation}
\min_{\bm{x}} \quad f\left(\bm{x}\right) + \mathbb{E}_{\bm{\xi}}\left[Q\left(\bm{x},\bm{\xi}\right)\right]\,,
\end{equation}
where $\bm{x}$ are the first stage decision variables, $Q\left(\bm{x},\bm{\xi}\right)$ is the second-stage cost function defined by,
\begin{align}
Q\left(\bm{x,\xi}\right) =  \min_{\bm{y}} \quad & q\left(\bm{x},\bm{\xi}\right)  \nonumber \\
\text{s. t.} \quad & \bm{T_{\xi}x + W_{\xi}y=b_{\xi}}, 
\end{align}
and $\bm{y}$ are the second-stage variables, and $\bm{\xi}$ are random variables.  General constraint equations are represented by $\bm{T_{\xi}}$, $\bm{W_{\xi}}$ and $\bm{b_{\xi}}$, and will be made specific to grid applications in the next section.  We restrict ourselves to using the linear stochastic programming problem, i.e.
\begin{equation}
\min_{\bm{x}} \quad \bm{c}^T\bm{x} + \mathbb{E}_{\bm{\xi}}\left[L\left(\bm{x},\bm{\xi}\right)\right]\,,
\end{equation}
where the recourse function is defined as the solution to,
\begin{align}
L\left(\bm{x,\xi}\right) =  \min_{\bm{y}} \quad & \bm{c_{\xi}}^{\sf T}\bm{y}  \nonumber \\
\text{s. t.} \quad & \bm{T_{\xi}x + W_{\xi}y=b_{\xi}}\\
& \bm{y}\geq \bm{0} \,. \nonumber
\end{align}
In the context of economic dispatch, the first stage variables $\bm{x}$ represent the dispatch levels of all generators, $\bm{c}$ are linear costs, and the expectation $\mathbb{E}_{\bm{\xi}}\left[\cdot\right]$ is taken with respect to the random variable $\bm{\xi}$, which represents errors in forecasting the available wind generation.

The second stage recourse variable, $\bm{y}$, represents adjustments to generation using reserves, such as up/down regulation, flex reserve or automatic generation control (AGC), that are necessary based on the amount of wind generation that materializes.

The optimization problem embedded inside the expectation $\mathbb{E}_{\bm{\xi}}\left[\cdot\right]$ makes it difficult to solve, so a sample average approximation (SAA) is commonly used to represent the expectation.  In the SAA approach, we introduce a set of $N$ scenarios, $\mathcal{S} = \left\lbrace\bm{\xi}_i\right\rbrace_{i=1}^N$, and approximate the expectation as $\mathbb{E}_{\bm{\xi}}\left[L\left(\bm{x},\bm{\xi}\right)\right] \approx \frac{1}{N}\sum_{i=1}^{N}L\left(\bm{x},\bm{\xi}_i\right)$.  Commuting the summation and minimization, we can simplify to the extensive form of a two-stage stochastic programming problem:
\begin{eqnarray}
\min_{\bm{x,y}_i i=1,\ldots,N} \quad \bm{c}^{\sf T}\bm{x} + \frac{1}{N}\sum_{i=1}^{N}\bm{c}_{\bm{\xi}_i}^{\sf T}\bm{y}_i \hspace{1.55in} & \nonumber \\
\begin{aligned}
\text{s. t.} \quad & \bm{Ax} & & & & & &=\bm{b} \\
& \bm{T}_{\bm{\xi}_1}\bm{x} \quad + & \bm{W}_{\bm{\xi}_1}\bm{y}_1 & & & & &= \bm{b}_1 \\
& \bm{T}_{\bm{\xi}_2}\bm{x} \quad + & & \bm{W}_{\bm{\xi}_2}\bm{y}_2  & & & &= \bm{b}_2 \\
& \bm{T}_{\bm{\xi}_3}\bm{x} \quad + & & & \ddots  & & &= \vdots \\
& \bm{T}_{\bm{\xi}_N}\bm{x} \quad + & & &  & \bm{W}_{\bm{\xi}_N}\bm{y}_N & &= \bm{b}_N
\\
& \bm{x}\geq\bm{0}, & \bm{y}_1 \geq \bm{0}, & \quad \bm{y}_2 \geq \bm{0}, & \ldots \:, & \bm{y}_N\geq\bm{0} \,.
\end{aligned}
\end{eqnarray}
	
A major challenge in the SAA approach is selecting the scenarios $\mathcal{S}$ such that the expectation $\mathbb{E}_{\bm{\xi}}\left[L\left(\bm{x},\bm{\xi}\right)\right]$ is accurately and efficiently calculated, while also ensuring reliable operation of the grid during different levels of wind generation.  In Monte Carlo approaches, the samples are drawn from their nominal distribution $p\left(\bm{\xi}\right)$ and given equal weights.  This approach yields $\mathcal{O}\left(N^{-1/2}\right)$ convergence, requiring prohibitively many samples to achieve an acceptable error level.  Other methods for selecting scenarios have been recently proposed, including scenarios motivated by grid reserve requirements \cite{papavasiliou2011reserve}, chance constraints \cite{cheng2018chance,lubin2016robust}, polynomial chaos expansions \cite{safta2017efficient}, or epi-spline basis functions \cite{staid2017generating}.

\subsection{Two-Stage Stochastic Economic Dispatch}
In this is study we restrict ourselves to studying the stochastic economic dispatch problem, i.e. balancing the (stochastic) demand for power across a grid with the (stochastic) power generation from all generators on the network. This study considers aggregated wind and load data, and therefore our model is formulated without network constraints. However, the scenario creation methods presented in this study should generalize to more complicated constraints on the network, e.g. DC and ACOPF. For a review of the economic dispatch problem, as well as the DCOPF and ACOPF problems, we recommend \cite{wood2012power,zhu2015optimization,taylor2015convex}. For an introduction to stochastic grid operations we recommend \cite{conejo2010decision}.
 
In the context of grid operations, the stochastic economic dispatch problem takes the following form:
\begin{align}
\min_{\substack{\bm{x} }} \quad & \sum_{\substack{g\in G}} \left( c_g x_{g} +\mathbb{E}_{\bm{\xi}}\left[L\left(\bm{x},\bm{\xi}\right)\right] \right) \nonumber  \\
\text{s. t.} \quad & x_g^{min} \leq x_{g} \leq x_g^{max} \quad \forall \: g \\
& -R_g^{down}\leq x_g - I_g \leq R_g^{up} \quad \forall \: g, \nonumber
\end{align}
where $\bm{x} = x_{g}$ for $g$ in the set of generators, $G$, are the first stage decisions representing thermal dispatch, and $L\left(\bm{x},\bm{\xi}\right)$ is the recourse function for a particular realization of the random variable $\bm{\xi}$, which represents the error in forecasting wind generation at each wind plant $w$ in the set of wind plants $W$.  The recourse function is defined as the solution to the following optimization problem:
\begin{align}\label{eq:loss-function}
&L\left(\bm{x,\xi}\right) =  \nonumber \\ 
\min_{\substack{\bm{y}^{\pm},\bm{\omega},\bm{\omega}^{spl} }}\quad  &\sum_{\substack{w\in W}} \left(c_w \omega_{w} + c_w^{spl} \omega_{w}^{spl} \right) \quad + \nonumber \\
& \sum_{g\in G} \left( c_q^+ y^+_{q} + c_q^- y^-_{q} \right) \quad + \nonumber \\
\text{s. t.} \quad  
& 0\leq y^+_{q} \quad \forall \: q \\
& 0\leq y^-_{q} \quad \forall \: q \nonumber \\
& 0  \leq  \omega_{w}  \leq \omega^{fcst}_{w} + \xi_{w}  \quad \forall \: w \nonumber \\
& \omega_{w}^{spl} = \left(\omega^{fcst}_{w} + \xi_{w}\right) \nonumber \\
& \qquad - \left( \omega_{w} \right)  \quad \forall \: w \nonumber \\
& \sum_{q \in D} \left(y^+_{q} - y^-_{q}\right) + \sum_{w \in W} \omega_{w} + \sum_{g\in G} x_{g}  \nonumber \\ 
& \qquad =\sum_{q \in D}\left(d_{q}  \right), \nonumber
\end{align}
where $\bm{\xi}=\xi_w, \, \bm{\omega}=\omega_w, \, \bm{\omega}^{spl}=\omega^{spl}_w \text{ for } w\in W$ and $\bm{y}^{\pm}=y^{\pm}_q$ for each bus $q$ in the set of buses $D$. The expectation $\mathbb{E}_{\bm{\xi}}\left[\cdot\right]$ is taken with respect to the random variable $\bm{\xi}$.  Generator ramp constraints are imposed between sequential timesteps, and not allowed to participate in recourse.  The loss of load or excess capacity represented by $\bm{y}^{\pm}$ is a slack variable that facilitates finding feasible solutions and represents a challenge to grid reliability and incurs substantial costs.

We assume the wind forecast errors are relative to a persistence forecast based on how well this approach performs in forecasting over short time-periods (see e.g. \cite{giebel2003state}). For an in-depth survey of wind power forecasting methods we recommend \cite{foley2012current} to the reader.

\section{Strategies for Scenario Generation}

\subsection{Importance Sampling}
Solving the two-stage stochastic economic dispatch problem requires a method of selecting scenarios. A logical starting point would be to draw from the nominal distribution of forecast errors in the wind power. However, this approach will not take into account costs communicated by the loss function \eqref{eq:loss-function}, potentially resulting in scenarios that do not include information about expensive events. As an example, consider a single wind farm generating power for one time period. Assuming that the nominal distribution of forecast errors is symmetric about zero, simply sampling from this distribution will produce scenarios that represent positive forecast errors, i.e. more wind than forecast, about half the time. If we use a small number of scenarios, e.g. $5$ to $10$, as is common in stochastic optimization problems, the scenario set may not include scenarios that represent negative forecast errors. This would be a costly mistake since both loss-of-load and reserves are more expensive than dispatching extra generation during the first stage decision process. 

To address this issue with random sampling we employ importance sampling, see e.g. \cite{bucklew_introduction_2004}. The reason for using importance sampling is that by including information from the cost function $L\left(\bm{x},\bm{\xi}\right)$, we construct a new probability density function that assign more mass to values of $\bm{\xi}$ corresponding to errors in overpredicting the available wind generation. Additionally, the construction of the importance distribution also provides a formula for properly weighting these samples to ensure the weighted average actually converges to $\mathbb{E}_{\bm{\xi}}\left[L\left(\bm{x},\bm{\xi}\right)\right]$.  Finally, importance sampling can also yield a lower variance estimate of the expected recourse cost.

Recall that our expected recourse cost is
\begin{equation}
\mu = \mathbb{E}_{\bm{\xi}}\left[L\left(\bm{x},\bm{\xi}\right)\right]=\int_{\Omega} L\left(\bm{x},\bm{\xi}\right)p\left(\bm{\xi}\right)\mathrm{d}\bm{\xi}
\end{equation}
where $p$ is a probability density function on $\mathbb{R}^d$ and $p\left(\bm{\xi}\right)=0$ for all $\bm{\xi}\notin \Omega$.  Introducing a new pdf $q$ on $\mathbb{R}^d$ (the importance distribution) we can write
\begin{align}
\mu = \mathbb{E}_{\bm{\xi}}\left[L\left(\bm{x},\bm{\xi}\right)\right] & =\int_{\Omega} L\left(\bm{x},\bm{\xi}\right)p\left(\bm{\xi}\right)\mathrm{d}\bm{\xi} \\
& = \int_{\Omega} \frac{L\left(\bm{x},\bm{\xi}\right)p\left(\bm{\xi}\right)}{q\left(\bm{\xi}\right)}q\left(\bm{\xi}\right)\mathrm{d}\bm{\xi} \\
& = \mathbb{E}_q\left[\frac{L\left(\bm{x},\bm{\xi}\right)p\left(\bm{\xi}\right)}{q\left(\bm{\xi}\right)}\right]
\end{align}
where $\mathbb{E}_q\left[\cdot\right]$ denotes expectation with $\bm{\xi}$ distributed as $q$.  The ratio $p\left(\bm{\xi}\right)/q\left(\bm{\xi}\right)$ is an adjusted weight on sample evaluations of $L$ that reflects sampling from $q$ instead of $p$.  

The importance sampling estimate of $\mu$ is then
\begin{equation}
\widehat{\mu}_q = \frac{1}{N}\sum_{i=1}^{N}\frac{L\left(\bm{x},\bm{\xi}_i\right)p\left(\bm{\xi}_i\right)}{q\left(\bm{\xi}_i\right)},
\end{equation}
where $\bm{\xi}_i$ are realizations of the random variable $\bm{\xi}$ with distribution $q\left(\bm{\xi}\right)$. The optimal importance sampling distribution, defined as $q\left(\bm{\xi}\right)=L\left(\bm{x},\bm{\xi}\right)p\left(\bm{\xi}\right)/\mu$, is a zero-variance estimator, and computes the correct expectation from a single sample.  Such an estimator suffers from the \textit{curse of circularity} (the quantity being estimated, $\mu$, is required to properly normalize the importance distribution), but provides guidance on constructing an approximate importance sampling distribution. 

Our approach to constructing an importance distribution is to solve a deterministic economic dispatch problem to approximate $\mu$ and $L\left(\bm{x},\bm{\xi}\right)$, and then use them to build an approximation of $q\left(\bm{\xi}\right)$. Our first simplification involves approximating the loss function as $L\left(\bm{x},0\right)$, i.e. we assume that the wind forecast is perfect. This step reduces the two stage stochastic optimization problem to a simple linear program:
\begin{align}
\bm{x}^* \: = \: \arg\min_{\substack{\bm{x} }} \quad & \sum_{\substack{g\in G}} \left( c_g x_{g} +L\left(\bm{x},\bm{0}\right) \right) \nonumber  \\
\text{s. t.} \quad & x_g^{min} \leq x_{g} \leq x_g^{max} \quad \forall \: g \\
& -R_g^{down}\leq x_g - I_g \leq R_g^{up} \quad \forall \: g, \nonumber
\end{align}
which we solve for $\bm{x}$. Given the solution $\bm{x}^*$ to the simplified dispatch problem we compute the expectation $\tilde{\mu}=\mathbb{E}\left[L\left(\bm{x}^*,\bm{\xi}\right)\right]$ via the trapezoidal rule. We remark that this process is fast since each evaluation of the loss function $L\left(\bm{x}^*,\bm{\xi}\right)$ is computed by solving a linear program. Given the approximated expected cost, $\tilde{\mu}$, the importance sampling distribution is defined as
\begin{equation}
q\left(\bm{\xi}\right)=\frac{L\left(\bm{x}^*,\bm{\xi}\right)}{\tilde{\mu}}p\left(\bm{\xi}\right)\,,
\end{equation}
and $q\left(\bm{\xi}\right)$ is then sampled instead of $p\left(\bm{\xi}\right)$ in the solution of the original stochastic economic dispatch problem.

\subsection{Gaussian Process Representation}
One shortcoming in the importance sampling approach presented in the preceding sub-section is that in reality, $\bm{x}$ is a first-stage variable in a two-stage stochastic optimization problem, and its value is computed at the same time as the solution to the rest of the problem. By ignoring this coupling, our proxy $\bm{x}^*$ may be off by a significant amount, leading to an importance distribution with an incorrect shape and location. 

We propose an alternative to using the importance sampling approach: representing the loss function as a Gaussian process (GP) and computing its expectation via a Bayesian quadrature technique. One benefit of this approach is that handles uncertainty in the second stage loss function introduced by not knowing the first stage dispatch.  In other words, we can treat $L\left(\bm{\xi}\right)$ as a random function when $\bm{x}$ is unknown.  Additionally, computing expectations of the GP with Bayesian quadrature does not require knowledge of the loss function to select the quadrature points. Instead, the Bayesian quadrature algorithm spreads out the sample points $\bm{\xi}_i$ in the domain based on the nominal distribution $p\left(\bm{\xi}\right)$ and the Gaussian process covariance kernel function.  The Bayesian quadrature algorithm also learns quadrature weights, thus ensuring that the weighted sum of the loss function sampled at the quadrature points returns the correct value for the expectation loss. A comprehensive review of Gaussian processes is provided by Rasmussen and Williams \cite{rasmussen_gaussian_2006}, while a description of the Bayesian quadrature sampling strategy can be found in \cite{ghahramani_bayesian_2003}.

By representing the loss function as a GP, we make the assumption that the joint distribution of any finite number of loss function evaluations is a multivariate Gaussian.  The GP is fully specified by a mean function $m\left(\bm{\xi}\right) = \mathbb{E}\left[L\left(\bm{\xi}\right)\right]$ and a kernel function $k\left(\bm{\xi},\bm{\xi}'\right) = \mathbb{E}\left[\left(L\left(\bm{\xi}\right)-m\left(\bm{\xi}\right)\right)\left(L\left(\bm{\xi^\prime}\right)-m\left(\bm{\xi^\prime}\right)\right)\right]$.  When represented as a GP, the loss function is
\begin{equation}
L\left(\bm{\xi}\right) \sim \mathcal{GP} \left[m\left(\bm{\xi}\right),k\left(\bm{\xi},\bm{\xi}'\right)\right]\,,
\end{equation}
and returns a Gaussian distribution when evaluated at a particular value $L\left(\bm{\xi}^*\right)$.  
The kernel function describes the covariance between any two points in the GP and encodes structural properties of the underlying random function like the smoothness, periodicity, and stationarity.  The kernel function is a nonlinear measure of similarity between inputs, with the dot product kernel leading to standard linear regression.  In this study we use a common nonlinear kernel function, the Gaussian or squared exponential function:
\begin{equation}
k\left(\bm{\xi},\bm{\xi}'\right) = \tau^2\exp\left(-\frac{|\bm{\xi}-\bm{\xi}'|^2}{2l^2}\right)
\end{equation}
where $\tau$ describes the magnitude and $l$ describes the smoothness of the function or correlation length between two points.

The GP can be conditioned on a set of loss function training points $\bm{\Xi}=\left\lbrace\bm{\xi}_i\right\rbrace_{i=1}^N$ to improve the posterior estimate of the loss function at new test points.  The joint distribution of training outputs $\mathbf{L}=L\left(\bm{\Xi}\right)$ and a test output $\mathbf{L}^*=L\left(\bm{\xi}^*\right)$ is
\begin{equation}
\left[\begin{array}{c} \mathbf{L} \\ \mathbf{L}^*\end{array}\right] \sim
\mathcal{N}\left(\left[\begin{array}{c} m\left(\bm{\Xi}\right) \\ m\left(\bm{\xi}^*\right)\end{array}\right], \left[\begin{array}{cc} k\left(\bm{\Xi},\bm{\Xi}\right) & k\left(\bm{\Xi},\bm{\xi}^*\right) \\ k\left(\bm{\xi}^*,\bm{\Xi}\right) & k\left(\bm{\xi}^*,\bm{\xi}^*\right)\end{array}\right] \right) \,.
\end{equation}
Conditioning on the training data $\mathcal{D}=\left \lbrace\bm{\Xi},\mathbf{L}\right\rbrace$ yields the following estimate for $\mathbf{L}^*$:
\begin{gather}
L\left(\bm{\xi}^* | \mathcal{D}\right)  \sim \mathcal{N}\left[m_{\mathcal{D}}\left(\bm{\xi}^*\right), \text{cov}_{\mathcal{D}}\left(L\left(\bm{\xi}^*\right),L\left(\bm{\xi}\right)\right) \right]\\
m_{\mathcal{D}}\left(\bm{\xi}^*\right) = \mathbf{k}_*^{\sf T}\mathbf{K}^{-1}\mathbf{L} \\
\text{cov}_{\mathcal{D}}\left(L\left(\bm{\xi}^*\right),L\left(\bm{\xi}\right)\right) = k\left(\bm{\xi}^*,\bm{\xi}^*\right) - \mathbf{k}_*^{\sf T}\mathbf{K}^{-1}\mathbf{k}_*
\end{gather}
where we have introduced the shorthand $\mathbf{k}_*=k_{*,i}=k\left(\bm{\Xi}_i,\bm{\xi}^*\right)$ and $\mathbf{K}=K_{ij}=k\left(\bm{\Xi}_i,\bm{\Xi}_j\right)$.  Note that the mean prediction $m_{\mathcal{D}}\left(\bm{\xi}^*\right)$ is simply a linear combination of the training points $\mathbf{L}$ and requires inverting a matrix that scales with the number of training points, N.  Additionally, the posterior covariance is given by the prior covariance less the information gained from the training points.

\subsection{Bayesian Quadrature}
Recall that in stochastic economic dispatch we are faced with calculating the expectation of our GP loss function $\mu = \int_{\Omega} L\left(\bm{\xi}\right)p\left(\bm{\xi}\right)\mathrm{d}\bm{\xi}$.  In Bayesian quadrature \cite{ghahramani_bayesian_2003,ohagan_bayeshermite_1991,hennig_probabilistic_2015} the integral $\mu$ is modeled as a Gaussian random variable because integration is a linear transformation that preserves Gaussianity and the integrand, $L\left(\bm{\xi}\right)$, is a Gaussian process.  The density function $p\left(\bm{\xi}\right)$ must also be a Gaussian, however this can be achieved by re-weighting via importance sampling.

Bayesian quadrature allows us to express the expectation and variance of  the loss analytically, and then devise a sampling strategy to optimally reduce this variance.  Following the derivation in \cite{ghahramani_bayesian_2003}, the loss averaged over possible functions is simply the expectation of the mean GP loss function:
\begin{align}
\mathbb{E}_{L|\mathcal{D}} \left[\mu\right] &= \int\int L\left(\bm{\xi}\right)p\left(\bm{\xi}\right)\mathrm{d}\bm{\xi} \; p\left(L|\mathcal{D}\right) \mathrm{d}L\\
& = \int \left[ \int L\left(\bm{\xi}\right) p\left(L|\mathcal{D}\right) \mathrm{d}L\right] p\left(\bm{\xi}\right)\mathrm{d}\bm{\xi} \\
& = \int m_{\mathcal{D}}\left(\bm{\xi}\right)p\left(\bm{\xi}\right)\mathrm{d}\bm{\xi}\\
& = \int k\left (\bm{\xi},\bm{\Xi}\right)\mathbf{K}^{-1}\mathbf{L} \, p\left(\bm{\xi}\right)\mathrm{d}\bm{\xi}\,.
\end{align}
Similarly, the expected variance of the estimate with respect to possible functions is 
\begin{align}
\mathbb{V}_{L|\mathcal{D}} \left[\mu\right] &= 
 \int\left[ \int L\left(\bm{\xi}\right)p\left(\bm{\xi}\right)\mathrm{d}\bm{\xi} - \right. \nonumber \\
& \ \ \  \left. \int m_{\mathcal{D}}\left(\bm{\xi}'\right)p\left(\bm{\xi}'\right)\mathrm{d}\bm{\xi}' \right]^2 p\left(L|\mathcal{D}\right) \mathrm{d}L \\
&=
\int\int\int\left[ L\left(\bm{\xi}\right) - m_{\mathcal{D}}\left(\bm{\xi}\right)\right] \left[ L\left(\bm{\xi}'\right) - \right. \nonumber\\
& \ \ \ \left. m_{\mathcal{D}}\left(\bm{\xi}'\right)\right]p\left(L|\mathcal{D}\right) \mathrm{d}L \; p\left(\bm{\xi}\right) p\left(\bm{\xi}'\right)\mathrm{d}\bm{\xi}\mathrm{d}\bm{\xi}'  \\
&=\int \int \text{cov}_{\mathcal{D}}\left(L\left(\bm{\xi}\right),L\left(\bm{\xi}'\right)\right)p\left(\bm{\xi}\right) p\left(\bm{\xi}'\right)\mathrm{d}\bm{\xi}\mathrm{d}\bm{\xi}'\\
&= \int\int \left[ k\left(\bm{\xi},\bm{\xi}'\right) - \right. \nonumber\\
& \ \ \ \left. k\left (\bm{\xi},\bm{\Xi}\right) \, \mathbf{K}^{-1} \, k\left(\bm{\Xi},\bm{\xi}\right) \right]p\left(\bm{\xi}\right) p\left(\bm{\xi}'\right)\mathrm{d}\bm{\xi}\mathrm{d}\bm{\xi}'\,.
\end{align}

A useful sampling strategy is to choose the set of samples that minimizes the expected variance:
\begin{equation}
\bm{\Xi}^* = \text{arg}\min_{\bm{\Xi}} \quad \mathbb{V}_{L|\mathcal{D}} \left[\mu\right] \,.
\end{equation}
This  can be compactly represented as 
\begin{equation}
\bm{\Xi}^* = \text{arg}\min_{\bm{\Xi}} \quad Z - \mathbf{w}^{\sf T} \mathbf{K}^{-1}\mathbf{w} \,,
\end{equation}
where we have defined the following quantities
\begin{gather}
\mathbf{w}=w_i := \int k\left(\bm{\xi},\bm{\Xi}_i\right)p\left(\bm{\xi}\right)\mathrm{d}\bm{\xi} \\
 Z = \int \int k\left(\bm{\xi},\bm{\xi}'\right) p\left(\bm{\xi}\right)p\left(\bm{\xi}'\right)\mathrm{d}\bm{\xi}\mathrm{d}\bm{\xi}'\,.
\end{gather}
Because the integrated prior covariance $Z$ has no dependence on the sample points, the optimization can be simplified to
\begin{equation} \label{eq:BQopt}
\bm{\Xi}^* = \text{arg}\max_{\bm{\Xi}} \quad\mathbf{w}^{\sf T} \mathbf{K}^{-1}\mathbf{w} \,,
\end{equation}
which shows that the optimal sampling strategy only depends on the choice of kernel covariance function $k\left(\bm{\xi},\bm{\xi}'\right)$ and the nominal distribution $p\left(\bm{\xi}\right)$, but not on the actual loss function evaluations $L\left(\bm{\Xi}\right)$.  To develop some intuition for Eq. \ref{eq:BQopt}, we can see that maximization involves a balance between selecting points that reflect the nominal distribution $p\left(\bm{\xi}\right)$ while avoiding redundancies from selecting similar points.

That the sampling strategy is independent of the loss function evaluations is a key differentiator between importance sampling and Bayesian quadrature approaches.  This is an advantage for Bayesian quadrature in cases where the loss function is unknown or stochastic, e.g. due to uncertainties in electrical loads.  Furthermore, as a quadrature rather than sampling strategy, the Bayesian quadrature points can be reused at successive timesteps if $k\left(\bm{\xi},\bm{\xi}'\right)$ and $p\left(\bm{\xi}\right)$ are unchanged.  The Bayesian quadrature estimate of the expected loss function takes the form of the following quadrature rule

\begin{equation}
\mathbb{E}_{L|\bm{\Xi},\bm{L}} \left[\mu\right] \approx \mathbf{w}^{\sf T}\mathbf{K}^{-1}\mathbf{L}=\sum^N_{i=1}\gamma_i L\left(\bm{\Xi}_i\right)
\end{equation}

\section{Performance of Different Sampling Strategies}

\subsubsection{Test data}

To study the performance of our economic dispatch methods we examine a flattened version of the Reliability Test System, modified by the Grid Modernization Lab Consortium (RTS-GMLC) \cite{GMLCRTS}. The RTS-GMLC system is an update to the original 1979 \cite{chairman1979ieee} and 1996 \cite{grigg1999ieee} version of the RTS that contains 73 buses, 96 conventional thermal generators and 4 wind plants. However, for our experiments we flatten the network and only consider power balance, i.e. matching the aggregated power output from the 96 conventional thermal generators and 4 wind plants to the aggregated demand without taking into consideration the network

The wind plant data used in our experiments is drawn from the Wind Integration National Dataset (WIND) Toolkit \cite{DRAXL2015355}. To appropriately model wind farms using WIND Toolkit data, wind plant power outputs were aggregated by first assigning geographic locations to nodes in the RTS-GMLC system and then summing the output power from the wind plants nearest to the node locations until the prescribed annualized capacity for the node in the RTS-GMLC system was reached.

The experiments in this section take place during a week in July. Load data from the RTS-GMLC system was combined with simulated wind data from 2013 from the WIND Toolkit. To determine $p\left(\bm{\xi}\right)$, we fit a Student's t-distribution to errors in a 5 minute persistence forecast.  A t-distribution was chosen because of the non-normal, heavy tailed behavior of changes in wind generation. The observation heavy-tailed distributions provide a good fit to wind data has been noted previously, see, e.g. \cite{hodge2011wind}.  The distributions of wind forecast errors could be further refined by conditioning on current power, time of day, season, etc, however this study used a fixed t-distribution that best matched errors aggregated over the full month of July.

\begin{figure}[!t]
\centering
\includegraphics[width=3.0in]{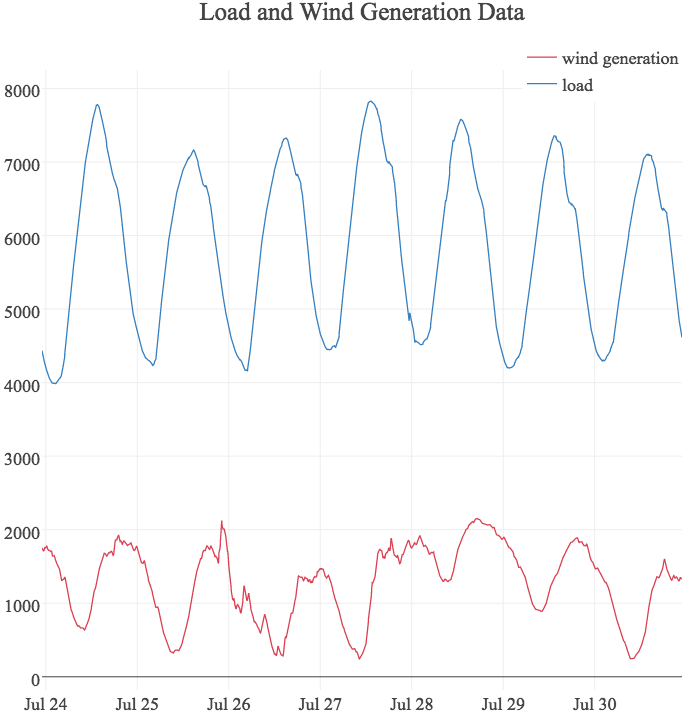}
\caption{Wind and load actuals used in economic dispatch experiments. We make note of the wind ramp event occurring on July 26th.}
\label{fig:load_and_generation}
\end{figure}

The actual loads and wind production data used in our experiments are shown in Figure~\ref{fig:load_and_generation}. There are two important aspects of the data: first, the peak renewable generation is just over one quarter of the peak demand, corresponding to a moderate renewable energy penetration scenario. Second, there is a wind ramp event that occurs around midnight July 25th. The wind ramp event will pose an economic and reliability challenge to test the different sampling strategies in our experiments.

\begin{figure}[!t]
\centering
\includegraphics[width=3.5in]{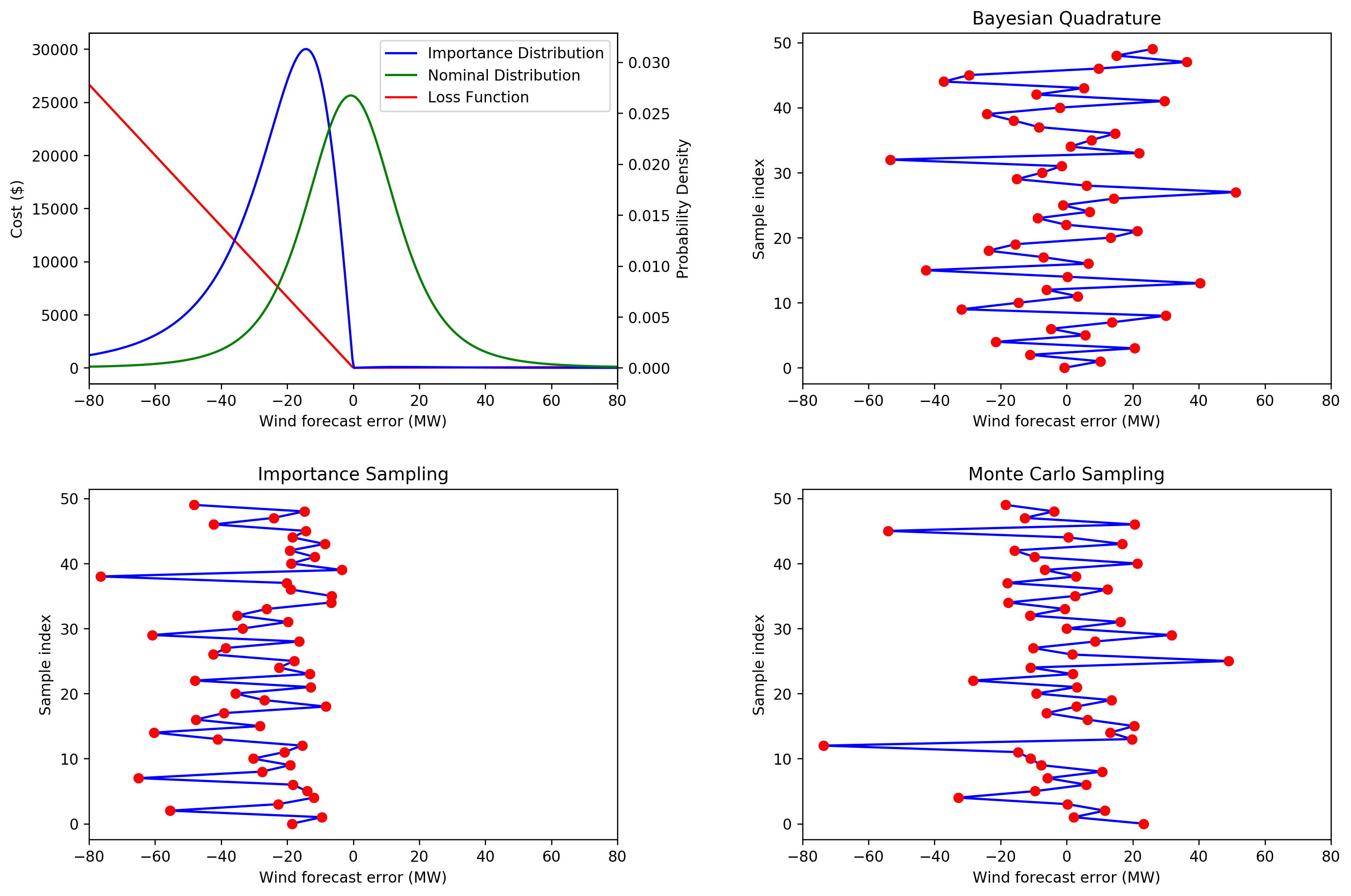}
\caption{The upper left panel shows the loss function (red), nominal distribution of $\xi$ (green), and importance distribution (blue) of forecast errors.  The other panels show representative samples generated by each sampling technique.}
\label{fig:importance-dist}
\end{figure}

\begin{figure*}[!t]
\centering
\includegraphics[width=0.9\textwidth]{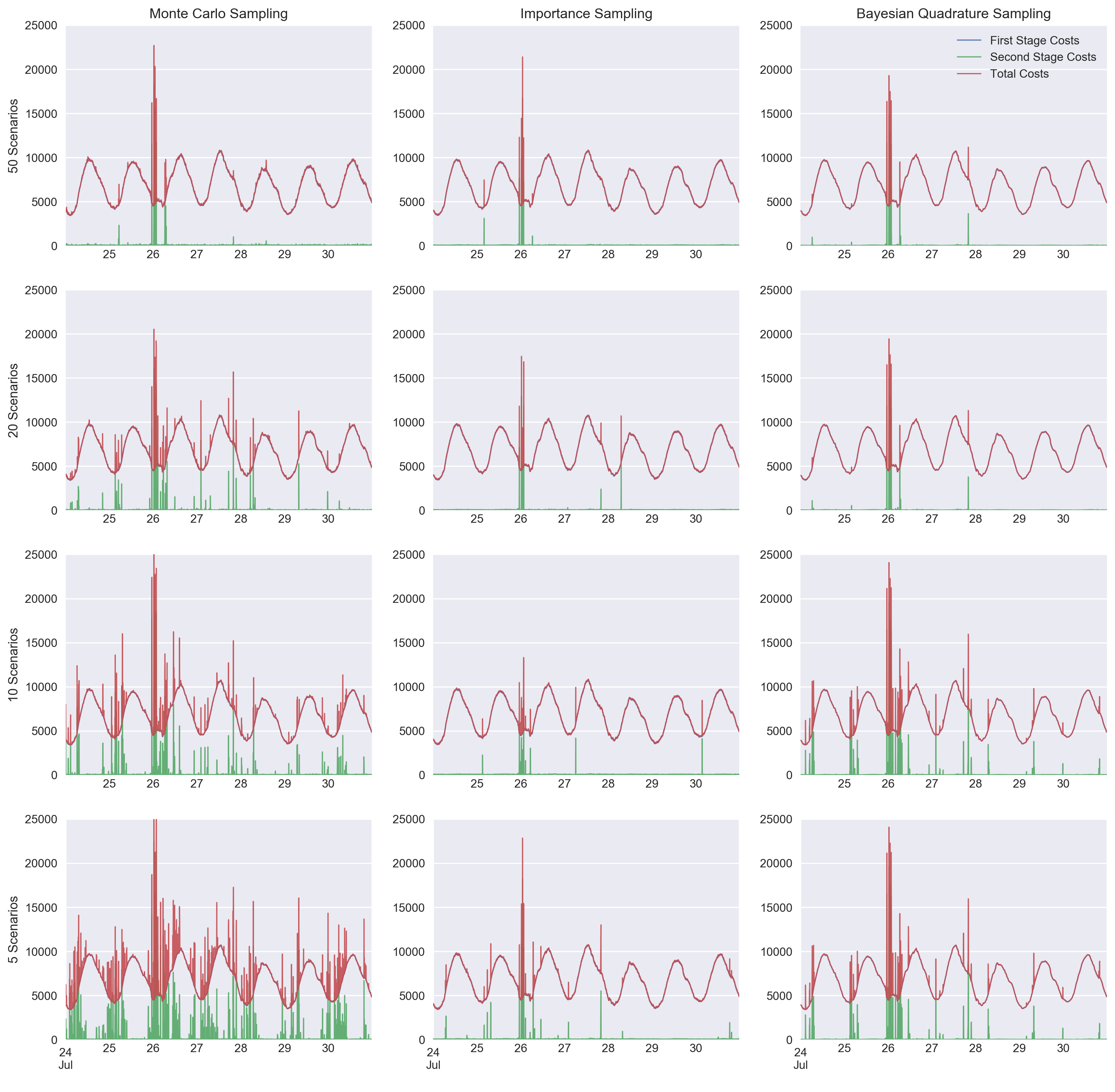}
\caption{Costs of the economic dispatch decisions from the three sampling strategies for 5, 10, 20, and 50 scenarios.  While the Bayesian quadrature approach triggers expensive second stage costs more often than importance sampling, it achieves lower total costs for 20 and 50 scenarios with much lower computational costs.}
\label{fig:scenarioTimeseries}
\end{figure*}

We examined the performance of each sampling strategy with 5, 10, 20, and 50 scenarios drawn at each 5 minute dispatch decision interval.  The resulting costs of the first and second stage decisions for the full week are shown in Figure ~\ref{fig:scenarioTimeseries} for all three strategies and the four sets of scenarios.  Additionally, a summary of the costs for the full week is shown in Table~\ref{tab:costs}.  To make the differences in sampling strategies more clear, our model does not explicitly implement reserves.  Instead any shortfall in generation results in an expensive loss-of-load penalty incurred during the second stage.  This results in the bulk of the costs attributed to first stage dispatch decisions, and rare, but comparatively large second stage costs.  In practice, these loss-of-load events would largely be addressed by ISO policies on reserve requirements.

Representative samples produced by each sampling technique, along with the nominal distribution, importance distribution, and loss function, are shown in Figure \ref{fig:importance-dist}.   The outcomes of these sampling strategies are shown in Figure \ref{fig:scenarioTimeseries} which displays a week-long timeseries of first and second stage dispatch cost.  The costs broken down by stage are additionally summarized in Table \ref{tab:costs2} and Table \ref{tab:costs3}.

The Monte Carlo method that draws from the nominal distribution has the most loss-of-load events. Alternatively, drawing from the importance sampling distribution leads to the smallest number of loss-of-load events. The reason for this behavior is clear: the shape of the loss function, large penalties for loss-of-load and relatively small penalties for overload or wind spilling, causes the importance distribution to favor sampling scenarios of wind overprediction. Thus, we are more likely to draw scenarios that represent lower than forecast values for wind from the importance distribution. These scenarios allow the stochastic optimization algorithm to hedge against expensive loss events. Drawing from the nominal distribution with Monte Carlo sampling, however, leads to scenarios that represent up and down ramps in wind power generation with equal probability. Therefore, for many of the 5 minute time-periods, the stochastic optimization algorithm was not equipped to hedge against loss-of-load events with Monte Carlo samples. 

The Bayesian quadrature approach results in a consistent decrease in maximum second stage costs as the number of scenarios is increased.  This consistent performance improvement is a strength of a quadrature-based approach.  In contrast, importance sampling does not show a monotonic decrease in peak second stage costs as more scenarios are considered.  This is likely due to errors arising from approximating the true importance distribution.

\begin{table}[!t]
\caption{Economic Dispatch Total Costs}
\label{tab:costs}
\centering
\begin{tabular}{cccc}
\hline
& \multicolumn{3}{c}{Sampling Strategy Costs (\$)} \\
\cline{2-4}
\# of Scenarios    & MC & IS & BQ\\
\hline
5      & \num{1.453e7}   & \num{1.383e7} & \num{1.394e7}      \\
10       & \num{1.409e7}   & \num{1.382e7} & \num{1.394e7}      \\
20       & \num{1.391e7}   & \num{1.383e7} & \num{1.378e7}      \\
50      & \num{1.386e7}   & \num{1.382e7} & \num{1.378e7}      \\
\hline
\end{tabular}
\end{table}

\begin{table}[!t]
\caption{Economic Dispatch First Stage Costs}
\label{tab:costs2}
\centering
\begin{tabular}{cccc}
\hline
& \multicolumn{3}{c}{Sampling Strategy Costs (\$)} \\
\cline{2-4}
\# of Scenarios    & MC & IS & BQ\\
\hline
5      & \num{1.355e7}   & \num{1.363e7} & \num{1.356e7}      \\
10       & \num{1.358e7}   & \num{1.365e7} & \num{1.356e7}      \\
20       & \num{1.360e7}   & \num{1.366e7} & \num{1.360e7}      \\
50      & \num{1.365e7}   & \num{1.366e7} & \num{1.360e7}      \\
\hline
\end{tabular}
\end{table}

\begin{table}[!t]
\caption{Economic Dispatch Second Stage Costs}
\label{tab:costs3}
\centering
\begin{tabular}{cccc}
\hline
& \multicolumn{3}{c}{Sampling Strategy Costs (\$)} \\
\cline{2-4}
\# of Scenarios    & MC & IS & BQ\\
\hline
5      & \num{9.839e5}   & \num{2.010e5} & \num{3.838e5}      \\
10       & \num{5.101e5}   & \num{1.688e5} & \num{3.838e5}      \\
20       & \num{3.073e5}   & \num{1.732e5} & \num{1.858e5}      \\
50      & \num{2.099e5}   & \num{1.648e5} & \num{1.834e5}      \\
\hline
\end{tabular}
\end{table}

With only 5 or 10 scenarios to consider, the importance sampling approach yields the lowest total costs.  We attribute this performance to the additional information available from the deterministic solution.  Recall that an estimate of the loss function $L\left(\bm{x},\bm{\xi}\right)$ is required to determine the importance distribution.  We approximate the loss function by solving a deterministic dispatch problem by assuming a perfect wind forecast.  We then use the trapezoidal rule to calculate the normalizing factor that makes the importance distribution a proper pdf.  Although this is a flawed estimate of the loss function, it is still valuable when only a few scenarios are available for optimization. 

For 20 or more scenarios, the Bayesian quadrature approach yields a lower total cost.  At this level of sampling, the Gaussian process approximation of the loss function is more accurate than the deterministic approximation made in importance sampling.  We also note that the Bayesian quadrature approach is much less computationally expensive than importance sampling because it avoids solving an additional deterministic problem and integrating to normalize the importance distribution.  Interestingly, the Bayesian quadrature solution yields a lower total cost than importance sampling despite triggering loss-of-load penalties more often.  This is a reminder that the economically optimal dispatch is not necessarily the most reliable.

\section{Conclusion}
This paper introduced two new techniques, importance sampling and Bayesian quadrature, to create scenarios of wind forecast error events for use in solving stochastic economic dispatch via a two-stage recourse problem.  The proposed techniques were compared against a Monte Carlo sampling strategy in a weeklong study of 5 minute economic dispatch using realistic generation, load, and wind forecast data from the RTS-GMLC test case and NREL's WIND Toolkit dataset.  Both importance sampling and Bayesian quadrature achieve a reduction in estimator variance as compared to the baseline Monte Carlo sampling, which ultimately yield lower total costs for the two-stage economic dispatch decisions.  Furthermore, the scenarios selected by the improved strategies are drawn further out into the tails of the distribution of possible wind ramp events than with Monte Carlo sampling.  The resulting dispatch decisions improve grid reliability and robustness by finding feasible solutions that reduce loss of load, while still maintaining economically optimal dispatch decisions.

We also found that the computational cost of determining sample locations is lower for Bayesian quadrature.  Importance sampling suffers from the \textit{curse of circularity} in that the normalization factor to make the importance distribution a proper probability density function requires knowing the expectation one is trying to estimate in the first place.  To overcome this, we solve a deterministic economic dispatch problem to estimate the loss function and then numerically integrate it to approximate the normalization constant.  This additional solve increases the computation cost of importance sampling relative to Bayesian quadrature.  Furthermore, in the Bayesian quadrature approach the sample points can be computed once and reused at successive timesteps because it is a true quadrature approach rather than a stochastic sampling approach.  The Bayesian quadrature samples only depend on the nominal wind pdf and the Gaussian process kernel, allowing them to be computed offline for wind distributions that are specific to a time of day, season, or current power output.  Our analysis simply used a Student's t-distribution to represent the distribution of wind ramps over the course of a year, however more condition-specific distributions would likely improve performance even further.

The most economic dispatch decisions are found when using importance sampling with ten or fewer scenarios, and with Bayesian quadrature for more than ten scenarios.  We attribute this to the extra information about the loss function that is gained through its deterministic approximation while computing the importance distribution.  With more than ten scenarios, the Bayesian quadrature approach can develop a better loss function surrogate and find a lower cost dispatch strategy.  While both proposed strategies improve the reliability from the baseline, the most economic strategy is not necessarily the most reliable.  The solutions found via Bayesian quadrature for 20 and 50 scenarios produce the lowest total cost as shown in Table I, but as shown in Table III, they involve more recourse or second stage costs than the less economic strategies found via importance sampling, reflecting increased number of expected loss of load events.  

Applying these strategies to a realistic test case of 5 minute economic dispatch with time varying loads, stochastic wind generation, and realistic costs yields improvements in total cost and in reliability represented by the second stage recourse costs.  Depending on the number of scenarios used, total costs are reduced by $4.8\%$ - $0.6\%$ relative to the Monte Carlo baseline.  A large component of the improvements in overall economic optimality is actually due to improvements in second stage costs that represent consumption of reserves or loss of load.  The economically optimal dispatch strategy reduces these costs by $79.6\%$ - $12.6\%$, demonstrating significant benefits to grid reliability and reserve requirements.  Employing these advanced scenario generation strategies can significantly reduce costs and improve the reliability of operating the grid with high levels of stochastic renewable energy production.

\section*{Acknowledgment}
This work was supported by the U.S. Department of Energy under Contract No. DE-AC36-08GO28308 with Alliance for Sustainable Energy, LLC, the Manager and Operator of the National Renewable Energy Laboratory. Funding provided by the Department of Energy's Grid Modernization Lab Consortium and Exascale Computing Project. 

\bibliography{bibliography}
\bibliographystyle{IEEEtran}

\end{document}